\documentclass[11pt]{amsart}
\usepackage{mathrsfs}

\normalparindent=\parindent

\theoremstyle{plain}
    \newtheorem{thm}{Theorem}
    \newtheorem{mainthm}[thm]{Main Theorem}
    \newtheorem{lem}[thm]{Lemma}
    
    \newtheorem{cor}[thm]{Corollary}
    \newtheorem{fact}[thm]{Fact}
\theoremstyle{definition}
    \newtheorem{defn}[thm]{Definition}
    
\theoremstyle{remark}
    \newtheorem{rem}[thm]{Remark}

\DeclareMathOperator{\mi}{min} \DeclareMathOperator{\ma}{max}
\DeclareMathOperator{\med}{med} \DeclareMathOperator{\maj}{maj}

\newcommand{\clf}[1]{\langle \{ #1 \}\rangle}

\newcommand{\C}{{\mathscr C}}
\newcommand{\F}{{\mathscr F}}
\newcommand{\J}{{\mathscr J}}

\newcommand{\OO}{{\mathscr O}}
\newcommand{\On}{{\mathscr O}^{(n)}}
\newcommand{\Ok}{{\mathscr O}^{(k)}}

\author[M.\,Pinsker]{Michael Pinsker}
\address{Department of Algebra and Computer Mathematics\\Vienna University of Technology}
\email{marula@gmx.at} \urladdr{http://www.algebra.tuwien.ac.at}
\title[Median functions]{The clone generated by the median functions}


\thanks{The author is supported by DOC [Doctoral Scholarship Programme of the Austrian Academy of Sciences].
He is grateful to M. Goldstern and L. Heindorf for their remarks
on the paper and to the II. Mathematisches Institut at Freie
Universit\"{a}t Berlin for their hospitality during his visit}

\begin{document}

\begin{abstract}
    Let $X$ be a linearly ordered set of arbitrary size (finite or infinite).
    Natural functions on such a set one can define using the linear order include maximum, minimum and median
    functions. While it is clear what the clone generated by the
    maximum or the minimum looks like, this is not obvious
    for median functions. We show that every clone on $X$ contains either no
    median function or all median functions, that is, the median functions generate each other.
\end{abstract}
\maketitle \thispagestyle{empty}

\begin{section}{Introduction}

\subsection{Preliminaries}
    Let $X$ be a set and denote by $\On$ the set of all $n$-ary functions
    on $X$. Then $\OO=\bigcup_{n=1}^{\infty}\On$ is the set of all functions
    on $X$. A \emph{clone} $\C$ over $X$ is a subset of $\OO$ which
    contains the projections and which is closed under compositions.
    The set of all clones over $X$ forms a complete lattice $Clone(X)$ with
    respect to inclusion. This lattice is a subset of the power set of $\OO$.

    A \emph{minimal} clone on $X$ is a minimal element in the lattice
    $Clone(X)\setminus \{\J\}$, where $\J$ is the trivial clone
    containing only the projections. Apparently every minimal
    clone is generated by a single nontrivial function. Functions which
    generate minimal clones are called \emph{minimal} as well.

\subsection{Notation}

    For a set of functions $\F$ we shall denote the smallest
    clone containing $\F$ by $\langle \F \rangle$.\\
    For a positive rational number $q$ we write
    $$
        \lfloor q \rfloor = \max\{n\in\mathbb{N}:\,n\leq q\}
    $$
    and
    $$
        \lceil q \rceil = \min\{n\in\mathbb{N}:\,q\leq n\}.
    $$
    If $a\in X^n$ is an $n$-tuple and $1\leq k\leq n$ we write $a_k$ for the $k$-th component of $a$.
    We will assume $X$ to be linearly ordered by the relation $<$. We let $\leq$ carry the obvious meaning.
\end{section}

\begin{section}{The median functions}
    Assume $X$ to be linearly ordered. We
    emphasize that the cardinality of $X$ is not relevant.
    For all $n\geq 1$ and all $1\leq k\leq n$ we define a function
    $$
        m^n_k(x_1,...,x_n)=x_{j_k}\quad ,\text{if} \,\, x_{j_1} \leq ... \leq
        x_{j_n}.
    $$
    In words, the function $m^n_k$ returns the $k$-th smallest element from an
    $n$-tuple. For example, $m^n_n$ is the maximum function $\max_n$
    and $m^n_1$ the minimum function $\min_n$ in $n$ variables. If
    $n$ is an odd number then we call $m^n_{\frac{n+1}{2}}$ the
    $n$-th median function and denote this function by $\med_n$.

    It is easy to check what the clones
    generated by the functions $\max$ and $\min$ look like:
    $$
        \clf{\ma_n}=\{\ma_k : k\geq 2\}\cup\J
    $$
    and
    $$
        \clf{\mi_n}=\{\mi_k : k\geq 2\}\cup\J
    $$
    where $n\geq 2$ is arbitrary. In particular, the two clones
    are minimal. Now it is natural to ask which of these properties hold for the
    functions ``in between'', that is the $m^n_k$ as defined before, most importantly the median
    functions.
    We will show that
    $$
        \clf{\med_n}\supseteq \{\med_k : k\geq 3 \text{ and
        odd}\}\cup\J
    $$
    but one readily constructs functions in that clone which are
    not a median function and not a projection. However, R. P\"{o}schel
    and L. Kalu\v{z}nin observed in \cite{PK79},
    Theorem 4.4.5, that the median of three variables (and hence by our result, all medians)
    does generate a minimal clone.
    \begin{fact}\label{medIsMinimal}
        The clone generated by the function $\med_3$ is minimal.
    \end{fact}
\subsection{Main theorem}
    \begin{mainthm}\label{medians}
        Let $k,n \geq 3$ be odd natural numbers. Then $\med_k\in \langle \{\med_n\} \rangle$. In
        other words, a clone contains either no median function or
        all median functions.
    \end{mainthm}

    We split the proof of the theorem into a sequence of
    lemmas.

    \begin{defn}
        Let $k,n\geq 1$ be natural numbers. Denote by $R(\frac{n}{k})$ the remainder of the
        division $\frac{n}{k}$. We say that $n$ is \emph{almost divisible by}
        $k$ iff either $R(\frac{n}{k})\leq \frac{n}{k}$ or $(k-1)-R(\frac{n}{k})\leq
        \frac{n}{k}$.
    \end{defn}

    Note that $n$ is almost divisible by $k$ if it is divisible by
    $k$. The following lemma tells us which medians of smaller
    arity are generated by $\med_n$ by simple identification of
    variables.

    \begin{lem}\label{almdivmed}
        Let $k\leq n$ be odd natural numbers. If $n$ is almost
        divisible by $k$, then $\med_k\in \langle \{\med_n\}\rangle$.
    \end{lem}

    \begin{proof}
        We claim that
        $$
            \med_k(x_1,...,x_k)=\med_n(x_1,...,x_1,x_2,...,x_2,...,x_k,...,x_k),
        $$
        where $x_j$ occurs in the $n$-tuple $\lfloor\frac{n}{k}\rfloor +1$ times if
        $j\leq R(\frac{n}{k})$ and $\lfloor\frac{n}{k}\rfloor$
        times otherwise. Assume $\med_k(x_1,...,x_k)=x_j$. Then
        there are at most $\frac{k-1}{2}$ components smaller than
        $x_j$ and at most $\frac{k-1}{2}$ components larger than
        $x_j$. Thus in our $n$-tuple, there are at most
        \begin{eqnarray}{\label{miss}}
        \frac{k-1}{2}\lfloor\frac{n}{k}\rfloor+\min(R(\frac{n}{k}),\frac{k-1}{2})
        \end{eqnarray}
        elements smaller (larger) than $x_j$.\\
        \textit{Case 1.}  $R(\frac{n}{k})\leq \frac{k-1}{2}$.\\ Since
        $n$ is almost divisible by $k$, we have either $R(\frac{n}{k})\leq
        \frac{n}{k}$ or $(k-1)-R(\frac{n}{k})\leq \frac{n}{k}$. In
        the latter case,
        $$
            R(\frac{n}{k})\leq \frac{k-1}{2} \quad \wedge \quad (k-1)-R(\frac{n}{k})\leq \frac{n}{k}
        $$
        and so
        $$
            R(\frac{n}{k})\leq \frac{n}{k}.
        $$
        Thus in either of the cases, we can calculate from
        (\ref{miss})
        \begin{eqnarray*}
        \begin{aligned}
            &\frac{k-1}{2}\lfloor\frac{n}{k}\rfloor+R(\frac{n}{k})\\
            &=\frac{1}{2}(k \lfloor\frac{n}{k}\rfloor +
            R(\frac{n}{k}))+\frac{1}{2}(R(\frac{n}{k})-\lfloor\frac{n}{k}\rfloor)\\
            &=\frac{n}{2}+\frac{1}{2}(R(\frac{n}{k})-\lfloor\frac{n}{k}\rfloor)\\
            &\leq\frac{n}{2}\\
        \end{aligned}
        \end{eqnarray*}
        and so $\med_n$ yields $x_j$.\\
        \textit{Case 2.} $\frac{k-1}{2}< R(\frac{n}{k})$.\\
        Again we know that either $R(\frac{n}{k})\leq
        \frac{n}{k}$ or $(k-1)-R(\frac{n}{k})\leq \frac{n}{k}$. In
        the first case, we see that
        $$
             \frac{k-1}{2}< R(\frac{n}{k}) \quad \wedge \quad R(\frac{n}{k})\leq \frac{n}{k}
        $$
        implies
        $$
            (k-1)-R(\frac{n}{k})\leq \frac{n}{k}
        $$
        and so (\ref{miss}) yields at most
        \begin{eqnarray*}
        \begin{aligned}
            &\frac{k-1}{2}\lfloor\frac{n}{k}\rfloor+\frac{k-1}{2}\\
            &=\frac{k-1}{2}\lfloor\frac{n}{k}\rfloor+\frac{1}{2}R(\frac{n}{k})+\frac{k-1}{2}-\frac{1}{2}R(\frac{n}{k})\\
            &=\frac{1}{2}(k \lfloor\frac{n}{k}\rfloor +
            R(\frac{n}{k}))-\frac{1}{2}\lfloor\frac{n}{k}\rfloor+\frac{k-1}{2}-\frac{1}{2}R(\frac{n}{k})\\
            &\leq\frac{n}{2}+\frac{1}{2}(-\lfloor\frac{n}{k}\rfloor+(k-1)-R(\frac{n}{k}))\\
            &\leq\frac{n}{2}\\
        \end{aligned}
        \end{eqnarray*}
        components which are smaller (larger) than $x_j$. This finishes the
        proof.
    \end{proof}
    \begin{cor}
        Let $k,n\geq 1$ be odd natural numbers. If $k\leq
        \sqrt{n}$, then $\med_k$ is generated by $\med_n$.
    \end{cor}
    \begin{proof}
        Trivially, $R(\frac{n}{k})\leq k-1$ and $k-1 \leq \frac{n}{k}$ as $k\leq
        \sqrt{n}$. Hence, $n$ is almost divisible by $k$.
    \end{proof}
    \begin{cor}\label{medn2med3}
        Let $n\geq 3$ be odd. Then $\med_3\in\clf{\med_n}$.
    \end{cor}
    \begin{proof}
        Simply observe that all $n\geq 4$ are almost divisible by $3$.
    \end{proof}
\subsection{Majority functions}
    We have seen that we can get small (that is, of small arity) median functions out of
    large ones. The other way we go with the help of majority
    functions.
    \begin{defn}
        Let $f\in\On$. We say that $f$ is a \emph{majority
        function} iff $f(x_1,...,x_n)=x$ whenever the value $x$
        occurs at least $\lceil\frac{n+1}{2}\rceil$ times among
        $(x_1,...,x_n)$.
    \end{defn}

    Note that $\med_n$ is a majority function for all odd $n$.
    We observe now that we can build a ternary majority function
    from most larger ones by identifying variables.

    \begin{lem}\label{LEM:ternaryMaj}
        Let $n\geq 5$ and let $\maj_n\in\On$ be a majority
        function. Then $\maj_n$ generates a majority function of
        three arguments.
    \end{lem}
    \begin{proof}
        Set
        $$
            \maj_3=\maj_n(x_1,...,x_1,x_2,...,x_2,x_3,...,x_3),
        $$
        where $x_j$ occurs in the $n$-tuple $\lfloor\frac{n}{3}\rfloor +1$ times if
        $j\leq R(\frac{n}{3})$ and $\lfloor\frac{n}{3}\rfloor$
        times otherwise. It is readily verified that $\maj_3$ is a
        majority function.
    \end{proof}

    The following lemma tells us that we can generate majority
    functions of even arity from majority functions of odd arity.

    \begin{lem}\label{lem:evenmaj}
        Let $n\geq 2$ be an even natural number. Then we can get
        an $n$-ary majority function $\maj_n$ out of any $n+1$-ary majority
        function $\maj_{n+1}$.
    \end{lem}
    \begin{proof}
        Set
        $$
            \maj_n(x_1,...,x_n)=\maj_{n+1}(x_1,...,x_n,x_n)
        $$
        and let $x\in X$ have a majority among $(x_1,...,x_n)$.
        Since $n$ is even, $x$ occurs $\frac{n}{2}+1$ times in the
        tuple which is enough for a majority in the $n+1$-tuple
        $(x_1,...,x_n,x_n)$.
    \end{proof}
    We now show that we can construct large majority functions out of
    small ones. This has already been known but we include an own proof
    here.
    \begin{lem} \label{maj}
        Let $n\geq 5$ be a natural number. Then we can
        construct an $n$-ary majority function out of any $(n-2)$-ary
        majority function $\maj_{n-2}$.
    \end{lem}
    \begin{proof}
        For $2\leq j\leq n-1$ and $1\leq i\leq n-1$ with $i\neq j$
        we define functions
        $$
            \gamma_i^j=
            \begin{cases}
                \maj_{n-2}(x_1,...,x_{i-1},x_{i+2},...,x_n)&,j\neq
                i+1\\
                \maj_{n-2}(x_1,...,x_{i-1},x_{i+1},x_{i+3},...,x_n)&,j=
                i+1\\
            \end{cases}
        $$
        In words, given an $n$-tuple $(x_1,...,x_n)$, $\gamma_i^j$ leaves away
        $x_i$ and the next component of the tuple which is not
        $x_j$ and calculates $\maj_{n-2}$ from what is left.
        Set
        $$
            z_j=\maj_{n-2}(\gamma_1^j,...,\gamma_{j-1}^j,\gamma_{j+1}^j,...,\gamma_{n-1}^j)
        $$
        and
        $$
            f=\maj_{n-2}(z_2,...,z_{n-1}).
        $$
        The function $f$ is an $n$-ary term of depth three over
        $\{\maj_{n-2}\}$.\\
        \textit{Claim.} $f$ is a majority function.\\
        We prove our claim for the case where $n$ is odd. The same proof works in the
        even case, the only difference being that the counting is slightly different (a majority occurs
        $\frac{n+2}{2}$ times instead of $\frac{n+1}{2}$ and so on). We leave the verification of this to the diligent
        reader.

        Assume $x\in X$ has a majority. If $x$ occurs more than
        $\frac{n+1}{2}$ times, then it is readily verified that all the $\gamma_i^j$ yield
        $x$ and so do all $z_j$ and so does $f$. So say $x$ appears
        exactly $\frac{n+1}{2}$ times among the variables of $f$.

        Next we observe that if $x_j=x$, then $z_j=x$: For if
        $\gamma_i^j\neq x$, then both components left away in $\gamma_i^j$, that is, $x_i$ and the component
        after $x_i$ which is not $x_j$, have to be equal to $x$. We can count
        $$
            |\{i:\gamma_i^j\neq x\}|\leq|\{i\neq j:x_i=x\}\setminus\{\max(i\neq j:
            x_i=x)\}|\leq\frac{n-1}{2}-1=\frac{n-3}{2}.
        $$
        Thus, $z_j=x$.

        Now we shall count a second time to see
        that if $x_1\neq x$ or $x_n\neq x$, then $f=x$: Say
        without loss of generality
        $x_1\neq x$. Then
        $$
            |\{2\leq j\leq n-1:
            x_j=x\}|\geq\frac{n+1}{2}-1=\frac{n-1}{2}
        $$
        and since we have seen that $z_j=x$ for all such $j$ we
        indeed obtain $f=x$.

        In a last step we consider the case where both $x_1=x$ and
        $x_n=x$. Let
        $$
            k=\min\{i: x_i\neq x\}
        $$
        and
        $$
            l=\max\{i:x_i\neq x\}.
        $$
        Since $n\geq 5$ those two indices are not equal. Count
        $$
            |\{i:\gamma^l_i\neq
            x\}|\leq|\{i:x_i=x\}\setminus\{k-1,n\}|=\frac{n+1}{2}-2=\frac{n-3}{2}.
        $$
        Thus, $z_l=x$ and we count for the last time
        $$
            |\{j: z_j=x\}|\geq|\{2\leq j\leq
            n-1:x_j=x\}\cup\{l\}|=\frac{n-3}{2}+1=\frac{n-1}{2}
        $$
        so that also in this case $f=x$.
    \end{proof}
    We conclude that if a clone contains a majority function, then
    it contains majority functions of all arities.
    \begin{cor}\label{allmaj}
        Let $n,k\geq 3$ be natural numbers. Assume $\maj_n\in\On$
        is any majority function. Then $\maj_n$ generates a majority function in $\Ok$.
    \end{cor}
    \begin{proof}
        If $k\geq n$ and $n,k$ are either both even or both odd, then we can iterate Lemma
        \ref{maj} to generate a majority function of arity $k$. Lemma \ref{lem:evenmaj}
        cares for the case when $k$ is even but $n$ is odd.

        In all other cases with $n\geq 5$, generate a ternary
        majority function from $\maj_n$ first with the help of Lemma
        \ref{LEM:ternaryMaj} and follow the procedure just
        described for the other case.

        Finally, if $n=4$, we can build a majority function $\maj_6$
        from $\maj_4$ first and are back in one of the other
        cases.
    \end{proof}
    Now we use the large majority functions to obtain large median functions.
    \begin{lem}
        For all odd $n\geq 3$ there exists $b\geq n$ such that
        $\med_n\in\clf{\med_3,\maj_b}$ for an arbitrary $b$-ary majority
        function $\maj_b$.
    \end{lem}
    \begin{proof}
        Let $n$ be given. Our strategy to calculate the median from an $n$-tuple will be the following: We apply
        $\med_3$ to all possible selections of three elements of the $n$-tuple. The results we write to an $n_1$-tuple,
        from which we again take all possible selections of three
        elements. We apply $\med_3$ again to these selections and
        so forth. Now the true median of the original $n$-tuple
        ``wins'' much more often in this procedure than the other
        elements, so that after a finite number of steps (a number we
        can give a bound for) more than half of the components of the then giant tuple
        have the true median as their value. To that tuple we
        apply a majority function and obtain the median.

        In detail, we define two sequences
        $(n_j)_{j\in\omega}$ and $(k_j)_{j\in\omega}$ by
        $$
            n_0=n,\quad n_{j+1}=\binom{n_j}{3}
        $$
        and
        $$
            k_0=1,\quad
            k_{j+1}=\binom{k_j}{3}+\binom{k_j}{2}(n_j-k_j)+\binom{k_j}{1}(\frac{n_j-k_j}{2})^2.
        $$
        The sequences have the following meaning: Given an
        $n_j$-tuple, there are $n_{j+1}$ possible selections of
        three elements of the tuple to which we apply the median $\med_3$. If the median of the $n_j$-tuple
        (which is equal to the median of the $n_0=n$-tuple) appeared $k_j$ times there, then it appears at least $k_{j+1}$
        times in the resulting $n_{j+1}$-tuple. Read $k_{j+1}$ as follows: We assume the worst case, namely that the median
        occurs only once in the original $n$-tuple, so $k_0=1$. If we pick
        three elements from the $n_j$-tuple and calculate $\med_3$, then the result is the median we are looking for
        if either all three elements are equal to the median ($\binom{k_j}{3}$ possibilities) or two are equal to the median
        ($\binom{k_j}{2}(n_j-k_j)$ possibilities) or one is equal to the median,
        one is smaller, and one is larger ($\binom{k_j}{1}\,(\frac{n_j-k_j}{2})^2$
        possibilities). Set
        $r_j=\frac{k_j}{n_j}$ for $j\geq 0$ to be the relative frequency of the median in the tuple after $j$ steps.
        We claim that $\limsup(r_j)_{j\in\omega}=1$:
        \begin{eqnarray*}
        \begin{aligned}
            r_{j+1}&=\frac{k_{j+1}}{n_{j+1}}\\
                   &=\frac{k_j}{n_j}\,\,\frac{(k_j-1)(k_j-2)+3(k_j-1)(n_j-k_j)+\frac{3}{2}(n_j^2-2n_jk_j+k_j^2)}{(n_j-1)(n_j-2)}\\
                   &=r_j\,\,\frac{3(n_j-1)^2+1-k_j^2}{2(n_j-1)(n_j-2)}
        \end{aligned}
        \end{eqnarray*}
        Further calculation yields
        \begin{eqnarray*}
        \begin{aligned}
            r_{j+1}&\geq r_j\,\,\frac{3(n_j-1)^2+1-k_j^2}{2(n_j-1)^2}\\
            &= r_j\,\,
            (\frac{3}{2}-\frac{(k_j-1)^2}{2(n_j-1)^2}-\frac{k_j-1}{(n_j-1)^2})\\
            &\geq r_j\,\,
            (\frac{3}{2}-\frac{1}{2}r_{j}^2-\frac{r_{j}}{n_j-1})\\
            &\geq r_j\,\,
            (\frac{3}{2}-\frac{1}{2}r_{j}^2-\frac{1}{n_j-1}).
        \end{aligned}
        \end{eqnarray*}
        Suppose towards contradiction that
        $(r_j)_{j\in\omega}$ is bounded away from $1$ by $p\,$: $r_j<p<1$ for all $j\in \omega$.
        Choose $j$ large enough so that
        $$\frac{1}{n_j-1}<\frac{1-p}{4}.$$
        Then
        \begin{eqnarray*}
        \begin{aligned}
            r_{i+1}&>r_i\,\,(\frac{3}{2}-\frac{p}{2}-\frac{1-p}{4})\\
            &=r_i\,\,(1+\frac{1-p}{4})
        \end{aligned}
        \end{eqnarray*}
        for all $i\geq j$ so that there exists $l>j$
        such that $r_l>p$, in contradiction to our assumption. Hence,
        $\limsup(r_j)_{j\in\omega}=1$.

        Now if we calculate $j$ such that $r_j>\frac{1}{2}$,
        we can obtain the median with the help of an $b=n_j$-ary majority function.
    \end{proof}
    We are ready to prove our main theorem.
    \begin{proof}[Proof of Main theorem \ref{medians}]
        Let $k, n$ be given. Corollary \ref{medn2med3} tells us that we can construct $\med_3$ out of $\med_n$.
        Since $\med_3$ is also a majority function, we can get majority
        functions of arbitrary arity with the help of Corollary \ref{allmaj}. Then by the preceding lemma, we can
        generate $\med_k$.
    \end{proof}
    \begin{rem}\label{REM:almostDivisibility}
        In fact, the lemma on almost divisibility is not needed
        for the proof of the theorem, since we only have to get
        $\med_3$ out of $\med_n$ (and $\med_3(x_1,x_2,x_3)=\med_n(x_1,x_2,...,x_2,x_3,...,x_3)$
        where $x_2$ and $x_3$ occur $\frac{n-1}{2}$ times
        in the $n$-tuple) and then apply Lemma \ref{maj} to generate large
        majority functions. Still, the lemma shows what we can
        construct by simple identification of variables.
    \end{rem}
\subsection{Minimality of the $m^n_k$}
    We mentioned that the clones generated by the maximum, the minimum and the median functions are
    minimal. Anyone who hoped that the same holds for all $m^n_k$ will be
    disappointed by the following lemma.
    \begin{lem}
        Let $n\geq 4$ and $2\leq k \leq \lfloor\frac{n}{2}\rfloor$.
        Then $m^n_k$ is not a minimal function.
    \end{lem}
    \begin{proof}
        It is enough to see that
        $$
            \mi_2(x,y)=m^n_k(x,...,x,y,...,y)\in\clf{m^n_k},
        $$
        where $x$ occurs in the $n$-tuple exactly
        $\lfloor\frac{n}{2}\rfloor$ times. The clone generated
        by $\min_2$ is obviously a nontrivial proper subclone of $\clf{m^n_k}$.
    \end{proof}
    Now comes just another disappointment.
    \begin{lem}
        Let $\lceil\frac{n}{2}\rceil< k < n$. Then $m^n_k$ is not
        a minimal function.
    \end{lem}
    \begin{proof}
        This time we have that
        $$
            \ma_2(x,y)=m^n_k(x,...,x,y,...,y)\in\clf{m^n_k},
        $$
        where $x$ occurs in the $n$-tuple exactly
        $\lfloor\frac{n}{2}\rfloor$ times. The clone generated by
        the maximum functions is obviously a proper subclone of
        $\clf{m^n_k}$.
    \end{proof}
    We summarize our results in the following corollary.
    \begin{cor}
        Let $n\geq 2$ and $1\leq k\leq n$. Then $m^n_k$ is minimal
        iff $k=1$ or $k=n$ or $n$ is odd and $k=\frac{n+1}{2}$.
        That is, the minimal functions among the $m^n_k$ are
        exactly the maximum, the minimum and the median functions.
    \end{cor}
\subsection{Remarks}
    For even natural numbers $n$ we did not define median
    functions. One could consider the so-called
    ``lower
    median'' instead:
    $$
        \med^{low}_n=m^n_{\frac{n}{2}}
    $$
    But as a consequence of the preceding corollary, $\med^{low}_n$ is
    not generated by the real medians and does therefore not serve
    as a perfect substitute. For the same reason, the ``upper
    median''
    $$
        \med^{upp}_n=m^n_{\frac{n}{2}+1}
    $$
    is not an ideal replacement either.

    However, the other direction almost works: $\med^{low}_n$
    generates the medians if and only if $n\geq 6$. Indeed,
    simple identification of variables suffices:
    $$
        \med_3=\med^{low}_n(x_1,...,x_1,x_2,...,x_2,x_3,...,x_3)
    $$
    where $x_j$ occurs in the $n$-tuple $\lfloor\frac{n}{3}\rfloor+1$ times
    if $j\leq R(\frac{n}{3})$ and $\lfloor\frac{n}{3}\rfloor$
    times otherwise. Of course we can do the same with the upper medians. It is easy to see that $\med^{low}_4$
    cannot generate the medians.

    One could have the idea of using a more general notion of
    median functions: Let $(X,\wedge,\vee)$ be a lattice. Define
    $$
        {\widetilde{m}}_k^n(x_1,...,x_n)=\bigwedge_{(j_1,...,j_k)\in N^k}\bigvee_{1\leq
        i\leq k}x_{j_i}.
    $$
    If the order induced by the lattice on $X$ is a chain,
    this definition agrees with our definition of $m^n_k$.
    However, although we can get $\widetilde{\med}_3$ out of $\widetilde{\med}_n$
    just like described in Remark \ref{REM:almostDivisibility}, our proof to
    obtain large medians via majority functions fails. We do
    not know under which conditions on the lattice the same results can be
    obtained.
\end{section}

\end{document}